\documentclass[12pt]{article}

\usepackage{amsmath,amssymb,theorem,makeidx,latexsym,epsfig}
\usepackage[normalem]{ulem}

\usepackage{cancel}

\newtheorem{defn}{Definition}[section]

\newtheorem{lemma}[defn]{Lemma}

{\theorembodyfont{\rmfamily}

\newtheorem{ex}[defn]{Example}}

\newtheorem{thm}[defn]{Theorem}

\newtheorem{prop}[defn]{Proposition}

\newtheorem{cor}[defn]{Corollary}

\newcommand{\h}{{\cal H}}

\newcommand{\ltp}{L^2(- \pi , \pi )}

\newcommand{\ltr}{ L^2(\mathbb R) }

\newcommand{\ltn}{{\ell}^2(\mathbb N)}

\newcommand{\si}{S^{-1}}

\newcommand{\afh}{ \forall f \in \h}

\newcommand{\mn}{\mathbb N}

\newcommand{\mr}{\mathbb R}

\newcommand{\mz}{\mathbb Z}

\newcommand{\mts}{ \{E_{mb}T_{na}g \}_{m,n \in \mz}}

\def\bp{{\noindent\bf Proof. \ }}

\def\ep{\hfill$\square$\par\bigskip}

\def\bqs{\begin{equation}}

\def\eqs{\tag*{$\square$}\end{equation}\par\bigskip}

\def\la{\langle}

\def\ra{\rangle}

\def\ga{\gamma}

\def\ftk{\{f_k\}_{k=1}^\infty}

\def\ctk{\{c_k\}_{k=1}^\infty}

\def\gtk{\left\{g_k\right\}_{k=1}^\infty}

\def\etk{\{e_k\}_{k=1}^\infty}

\def\suk{\sum_{k=1}^\infty}

\def\nl{\left|\left|}

\def\nr{\right|\right|}

\def\span{\overline{\text{span}}}

\def\supp{\text{supp}}

\def\bop{\begin{op}\rm}

\def\eop{\end{op}}

\def\bee{\begin{eqnarray}}

\def\ene{\end{eqnarray}}

\def\bes{\begin{eqnarray*}}

\def\ens{\end{eqnarray*}}

\def\bei{\begin{itemize}}

\def\eni{\end{itemize}}

\def\bt{\begin{thm}}

\def\et{\end{thm}}

\def\bc{\begin{cor}}

\def\ec{\end{cor}}

\def\bpr{\begin{prop}}

\def\epr{\end{prop}}

\def\bl{\begin{lemma}}

\def\el{\end{lemma}}

\def\bd{\begin{defn}}

\def\ed{\end{defn}}

\def\bex{\begin{ex}}

\def\enx{\end{ex}}

\def\bfi{\begin{fig}}

\def\efi{\end{fig}}

\def\inr{\int_{-\infty}^\infty}

\newcommand{\nft}{ || f||^2}

\def\hpi{\hat{\psi}}

\def\mzd{{\mathbb Z}^d}

\newcommand{\ltrd}{ L^2({\mathbb R}^d) }

\def\dy{\{D^jT_k\psi\}_{j,k\in \mz}}

\def\inr{\int_{-\infty}^\infty}

\def\nti{\{n_i\}_{i\in I}}

\def\eti{\{e_i\}_{i\in I}}

\def\hti{\{h_i\}_{i\in I}}

\def\fti{\{f_i\}_{i\in I}}

\def\gti{\{g_i\}_{i\in I}}

\def\oti{\{\omega_i\}_{i\in I}}

\def\otj{\{\omega_j\}_{j\in I}}

\def\ottk{\widetilde{\omega_k}}

\def\oj{\omega_j}

\def\sui{\sum_{i\in I}}

\def\suj{\sum_{j\in I}}

\def\wp{W^\bot}

\def\mtho{\{E_{mb}T_{na}h_1\}_{m,n\in \mz}}
\def\mtgo{\{E_{mb}T_{na}g_1\}_{m,n\in \mz}}
\def\mtht{\{E_{mb}T_{na}h_2\}_{m,n\in \mz}}
\def\mtro{\{E_{mb}T_{na}r_1\}_{m,n\in \mz}}
\def\mtrt{\{E_{mb}T_{na}r_2\}_{m,n\in \mz}}
\def\mtgt{\{E_{mb}T_{na}g_2 \}_{m,n\in \mz}}

\def\wpo{\{D^jT_k\psi_1\}_{j,k\in \mz}}
\def\wpt{\{D^jT_k\psi_2\}_{j,k\in \mz}}
\def\wpto{\{D^jT_k\widetilde{\psi_1}\}_{j,k\in \mz}}
\def\wptt{\{D^jT_k\widetilde{\psi_2}\}_{j,k\in \mz}}

\def\atj{\{a_j\}_{i\in J}}
\def\btj{\{b_j\}_{i\in J}}
\def\ptj{\{p_j\}_{i\in J}}
\def\qtj{\{q_j\}_{i\in J}}

\def\sumn{\sum_{m,n\in \mz}}

\def\hpi{\hat{\psi}}
\def\htpi{\hat{\tilde{\psi}}}
\def\tpi{\tilde{\psi}}

\title{\sout{Seven} Six problems in frame theory}

\date{\today}

\author{Ole Christensen }

\begin{document}

\maketitle


\noindent {\it Dedicated to Professor Butzer on the occasion of his eighty-fifth
anniversary} \\

\section{Introduction}

Progress in research is due to great scientists. It is due to great mentors and teachers who inspire and challenge their colleagues and students to consider new problems. And it is due to great personalities, whose positive attitude and support make us go the extra mile in order to reach the goal. Professor Butzer has all these properties. It is a great honor for me to be invited to contribute to a book that celebrates him and his work, and I want to thank Professor Schmeisser and Professor Zayed for giving me this opportunity. Finally, at last, but not least: Paul - thanks for all!

The purpose of this paper is to present selected problems in frame theory that have been open for some years.
We will also discuss the role of these problems and the technical difficulties that have prevented them from being solved so far. The hope is that the paper will contribute to the solution of at least some of these problems.

The paper is organized as follows. In Section \ref{404b} we provide a very short presentation of the necessary background from frame theory, operator theory, wavelet theory, and Gabor analysis. We only cover the parts that are essential in order to understand the open problems and their role in frame theory. The first of these is considered in Section \ref{404a}, where we deal with the issue of extending  pairs of Bessel sequences to pairs of dual frames. In
the abstract setting this can always be done, but when we ask for the extended systems to have a particular structure, various open problems appear.

In Section \ref{404c} we consider the theory for R-duals, introduced by Casazza, Kutyniok, and Lammers. Following their paper \cite{CKL} we discuss the question whether this theory leads to a generalization of the duality principle in Gabor analysis. We also describe an alternative approach that has been developed in \cite{CKK2}. Section \ref{404d} deals with the construction of wave packet frames in $\ltr.$ Such systems appear by the combined action of a class of translation, modulation, and scaling operators, and it turns out that the parameters in these operators have to be chosen carefully in order not to violate the Bessel condition. It is not clear yet how one can describe suitable choices in a general way.
In Section \ref{404e} we consider a concrete class of functions, namely, the B-splines, and raise the question of  the exact range of parameters for which they generate Gabor frames. This type of question has been considered in the literature for various other functions than the $B$-splines, but exact answers are only known  for relatively few functions. In Section \ref{404g}, two open problems concerning finite structured frames are considered: first, the famous Heil--Ramanathan--Topivala conjecture, stating that any finite nontrivial Gabor system is linearly independent, and secondly, the problem of finding good estimates for the lower frame bound for a finite collection of exponentials in $\ltp.$
Finally, Section \ref{404h} describes the Feichtinger conjecture (2003), which has attracted considerably attention over the last years. It is known to be equivalent to the Kadison--Singer problem dating back to 1959. The Feichtinger conjecture was reported to be solved affirmatively shortly before submission of this manuscript.

\section{Preliminaries} \label{404b}
The purpose of this section is to give a short presentation of frame theory, with focus on the parts of the theory that are necessary in order to understand the open problems. It is not the intention to give a full survey on frame theory. We refer to the books \cite{Y}, \cite{Da2}, and \cite{CBN} for more information.

\subsection{Frames and Riesz bases in Hilbert spaces}

In the entire section $(\h, \la \cdot, \cdot \ra)$ denotes a separable Hilbert space.

\newpage

 \bd A sequence $\ftk$ of elements in $\h$ is called a
Riesz sequence if there exist constants $A,B>0$ such that
\bes A\sum |c_k|^2\leq\nl \sum c_kf_k\nr^2\leq B\sum |c_k|^2
\ens for all finite sequences $\{c_k\}.$ A Riesz sequence $\ftk$ is
called a Riesz basis if
$\span \ftk=\h.$  \ed

Alternatively, Riesz bases are precisely the sequences that
have the form
$\ftk= \{Ue_k\}_{k=1}^\infty$ for some orthonormal basis  $\etk$ and some
bounded bijective operator $U:\h\to \h.$
For a given Riesz basis $\ftk= \{Ue_k\}_{k=1}^\infty,$ the sequence
$\gtk=\{(U^{-1})^*e_k\}_{k=1}^\infty$
is called the {\it dual} of  $\ftk= \{Ue_k\}_{k=1}^\infty.$ Note that the dual of $\gtk$ is \bes
\left\{\left(\left((U^{-1})^*\right)^{-1}\right)^*e_k\right\}_{k=1}^\infty
=\{Ue_k\}_{k=1}^\infty = \ftk,\ens i.e.,  $\ftk$ and $\gtk$ are duals of each
other. A Riesz basis $\ftk$ and its dual $\gtk$ are biorthogonal, i.e.,
\index{biorthogonal sequences}
\bes \la f_j,g_k\ra = \delta_{k,j}, \ j,k\in \mn.\ens

\bd A sequence $\ftk$ in $\h$ is called a Bessel
sequence \index{Bessel sequence}   if there exists a constant
$B>0$ such that \bes \suk |\la f,f_k\ra  |^2 \le B \ \nft , \
\afh. \ens \ed

Associated to a Bessel sequence
$\ftk,$ the {\it pre-frame operator} or {\it synthesis operator} is
\bes
T: \ltn \to \h , \ T\ctk = \suk c_kf_k. \ens The operator $T$ is bounded for any Bessel sequence $\ftk.$
The adjoint operator of $T$ is called the {\it analysis
operator} and is given by  \bes T^*: \h \to \ltn , \ \ T^* f= \{\la f,f_k \ra
\}_{k =1}^{\infty}.  \ens
Finally, the {\it frame operator} is defined by \index{frame
operator} \bes  S : \h \to \h , \ \ Sf=TT^*f= \suk \la f,f_k\ra
f_k . \ens

The classical definition of a frame, originally given by Duffin and Schaeffer
\cite{DS}, reads as follows.

\bd A sequence $\ftk$ in $\h$ is a {\it frame} if there exist
constants $A,B>0$ such that \bes A\,||f||^2 \le
\suk | \la f,f_k\ra|^2\le B \, ||f||^2, \ \forall f\in \h.\ens
$A$ and $B$ are called {\it frame bounds.}
A frame is {\it tight} if we can take $A=B.$ \ed

Note that any Riesz basis is a frame. A frame which is not a Riesz basis, is said to be {\it
overcomplete} or {\it redundant.} The following classical result  shows that any
frame leads to an expansion of the elements in $\h$ as a (infinite) linear combinations of the frame elements. It also shows that the general expansion simplifies considerably for tight frames. Finally, the last part of the result shows that for frames that are not Riesz bases, the coefficients in the series expansion of an element $f\in \h$ are not unique:

\bt
Let $\ftk$ be a frame with frame operator $S.$  Then the following hold:
\bei \item[(i)] Each $f\in \h$ has the decompositions
\bes f=  \suk \la  f, \si f_k \ra f_k=  \suk \la  f,  f_k \ra  \si f_k.\ens
\item[(ii)] If $\ftk$ is a tight frame with frame bound $A$, then $S=AI,$ and \bee
\label{7221b} f= \frac1{A} \suk \la f,f_k\ra f_k, \ \forall f\in
\h.\ene \item[(iii)] If $\ftk$ is an overcomplete frame.
there exist frames \\ $\gtk \neq \{\si f_k\}_{k=1}^\infty$ for
which \bee \label{173a}  f=\sum_{k=1}^\infty \la f, g_k\ra f_k, \ \forall f\in
\h.\ene \eni \et

Any Bessel sequence $\gtk$ satisfying \eqref{173a} for a given frame $\ftk$
is called a {\it dual frame} of $\ftk$. The special choice
$\gtk = \{ \si f_k\}_{k=1}^\infty$ is called the {\it canonical dual frame.} In order to avoid confusion we note that if \eqref{173a} holds for two Bessel sequences
$\ftk$ and $\gtk,$ they are automatically frames:

\bl \label{223a} If $\ftk$ and $\gtk$ are Bessel sequences and \eqref{173a} holds, then
$\ftk$ and $\gtk$ are dual frames. \el

\subsection{Operators on $\ltr$}

Most of the open problems considered in this paper will deal with Gabor systems or wavelets systems. In order to describe these systems we need the following unitary operators on $\ltr:$


\vspace{.1in}\noindent{\it Translation by} \ $a\in \mr$:
\  \  $T_a: \ltr \to \ltr, \
(T_af)(x)=f(x-a).$

\vspace{.3cm}\noindent{\it Modulation by} \ $b\in \mr:
\ \  E_b: \ltr \to \ltr, \
(E_bf)(x)=e^{2\pi ibx}f(x).$

\vspace{.3cm}\noindent{\it Dilation by} \ $a > 0: \ \hspace{.5cm} \
D_a: \ltr \to \ltr, \
(D_af)(x)= \frac{1}{\sqrt{a}}f(\frac{x}{a}).$

\vspace{.3cm}\noindent{\it Dyadic scaling:} \  \hspace{1cm} \
$D: \ltr \to \ltr, \
(Df)(x)=  2^{1/2}f(2x).$

\vspace{.2cm}\noindent All these operators are unitary on $\ltr.$ We will
also consider the {\it Fourier transform,} \index{Fourier
transform}  for $f \in L^1(\mr)$ defined by
$$ {\cal F}f(\ga)= \hat{f}(\ga):= \inr f(x)e^{-2\pi i x\ga}\,dx, \ \ga\in \mr,$$ and extended to a unitary operator on
$\ltr$ in the usual way. The operators $T_a, E_b, D,$ and
${\cal F}$ are related by the following commutator relations:
\bes T_aE_b & = & e^{-2\pi i ba}E_bT_a, \ \ T_bD= DT_{b/a}, \
\ DE_b= E_{b/a}D  \\
 {\cal F}T_a & = & E_{-a}{\cal F}, \ \ {\cal F}E_a=T_a{\cal
F}, \  \ {\cal
F}D=D^{-1}{\cal F}. \ens

\subsection{Wavelet systems}

A systems of functions on the form $\dy$ for a fixed function $\psi\in \ltr$ is called a {\it dyadic wavelet system.} Note that
\bes
D^jT_k\psi(x)=2^{j/2}\psi(2^jx-k), \ x\in \mr. \ens
Given a  frame $\dy$ for $\ltr,$
the associated frame
operator is \bes S: \ltr \to \ltr, \  Sf= \sum_{j,k\in
\mz} \la f, D^jT_k\psi \ra D^jT_k\psi,\ens and
the frame
decomposition reads \bes f=
\sum_{j,k\in \mz} \la f, \si D^jT_k\psi \ra D^jT_k\psi, \ f\in
\ltr.\ens
In order to use the frame decomposition we needs to calculate
the numbers $ \la f, \si
D^jT_k\psi \ra$ for all  $ j,k\in \mz,$ i.e., a double-infinite sequence of numbers. One can show that \bes \si D^jT_k\psi= D^j \si T_k\psi, \ens so in practice it is enough to calculate the action of $\si$ on the functions $T_k\psi,$ and then apply the scaling $D^j.$
Unfortunately, in general
\bes D^j \si T_k\psi \neq D^j T_k\si\psi.\ens

Thus, we
can not expect the canonical dual frame of a wavelet frame
to have wavelet structure. As a concrete example
(taken from \cite{Da1} and \cite{Chui}), let $\{D^jT_k\psi\}_{j,k\in \mz}$ be a wavelet
orthonormal basis for $\ltr$. Given $\epsilon \in ]0,1[$, let
$\theta = \psi + \epsilon D\psi.$
Then $\{D^jT_k\theta\}_{j,k\in \mz}$ is a
Riesz basis, but the canonical dual frame of
$\{D^jT_k\theta\}_{j,k\in \mz}$ does {\it not} have the wavelet
structure. Since the dual is unique for a Riesz basis,
this example demonstrates that there are wavelet frames where no dual with wavelet
structure exists. On the other hand Bownik and Weber  \cite{BW} have given
an interesting example of a wavelet frame
$\dy$ for which the canonical dual does
not have the wavelet structure, but other dual frames with wavelet structure exist. Dual pairs of wavelet frames can be characterized as follows, see \cite{CS}; we restrict us to the case $a=2, b=1,$ although a similar result holds for general
parameters.

\bt Two Bessel sequences
$\{D^jT_{k}\psi\}_{j,k\in \mz}$ and
$\{D^jT_{k}\widetilde{\psi}\}_{j,k\in \mz}$ form dual wavelet frames for $\ltr$ if and only if the following two conditions are satisfied:
\bei \item[(i)] $\sum_{j\in \mz}\overline{\widehat{\psi}(2^j\ga)}\widehat{\widetilde{\psi}}(2^j\ga)
=b$ for a.e.\ $\ga\in \mr.$
\item[(ii)] For any number $\alpha\neq 0$ of the
form $\alpha=m/2^j, \ m,j\in \mz,$ \bes
\sum_{\{(j,m)\in \mz^2 \: | \: \alpha =m/2^j\}} \overline{\widehat{\psi}(2^j\ga)}\widehat{\widetilde{\psi}}(2^j\ga+m)=0, \ a.e. \ \ga \in \mr.\ens
\eni \et

\subsection{Gabor systems} \label{273a}

Gabor systems in $\ltr$  have the form \bes \{e^{2\pi
imbx}g(x-na)\}_{m,n \in \mz}\ens for some $g\in \ltr, \ a,b>0$. Using operator
notation, we can write a Gabor system as
$\mts.$

We will not go into a general description of Gabor analysis and its role in time-frequency analysis, but just refer to the books \cite{FS1, FS2, G2}.

Letting $\chi_{[0,1]}$ denote the characteristic function for the interval $[0,1],$
it is easy to show that $\{E_mT_n \chi_{[0,1]} \}_{m,n\in \mz}$ is an orthonormal basis for $\ltr.$ But the function $\chi_{[0,1]}$ is discontinuous and
has very slow decay in the Fourier domain, so the function is not suitable for time-frequency analysis.  For the sake of time-frequency analysis we want the Gabor system $\mts$ to be
generated by a continuous function $g$ with compact support. The following classical result shows that this more or less forces us to work with frames.

\bl If $g$ is be a continuous function with compact support, then
\bei
\item[(i)]  $\{E_{mb}T_{na}g\}_{m,n\in \mz}$ can not be an ONB.

\item[(ii)]  $\{E_{mb}T_{na}g\}_{m,n\in \mz}$ can not be a Riesz basis.

\item[(iii)] $\{E_{mb}T_{na}g\}_{m,n\in \mz}$ can  be a frame if
$0<ab< 1;$
\eni\el

In addition to (iii), if $0< ab<1,$ it is always possible to find a function
$g\in C_c(\mr)$ such that $\mts$ is a Gabor frame. We also note that
no matter whether $g$ is continuous or not, Gabor frames $\mts$ for $\ltr$ only exist if $ab\le 1.$

Bessel sequences of the form $\mts$ will play a central role in some of the open problems to be considered in this article, so let us state a classical sufficient condition that is easy to verify.

\bl \label{804a} Let $g$ be a bounded function with compact support. Then $\mts$ is
a Bessel sequence for any  $a,b>0$.\el

We have already seen that for a wavelet system, the frame operator in general does not commute with the translation operator. For Gabor systems, the frame operator commutes with the relevant operators.
We will need the result below, which is almost identical to
Lemma 9.3.1 in \cite{CBN}.

\bl \label{223b} Let $g,h\in \ltr$ and $a,b>0$ be
given, and assume that $\mts$ and
$\{E_{mb}T_{na}h\}_{m,n\in \mz}$ are Bessel sequences. Then the following holds: \bei \item[(i)] Letting $T$ and $U$ denote the preframe operators for  $\mts,$ respectively
$\{E_{mb}T_{na}h\}_{m,n\in \mz},$
\bes TUE_{mb}T_{na}=E_{mb}T_{na} TU, \ \forall m,n \in \mz.\ens
\item[(ii)] If $\mts$ is a frame with frame operator
$S=TT^*,$ then also \bes
\si E_{mb}T_{na}=E_{mb}T_{na}\si, \ \forall m,n \in \mz.\ens
\eni  \el

Lemma \ref{223b}(ii) implies that for a Gabor frame $\mts$ with
associated frame operator $S$, the canonical dual frame also has Gabor structure, in contrast with the situation we encountered for wavelet frames. However, even for a nice frame $\mts$ it is nontrivial to control the properties of the canonical dual frame
$\{E_{mb}T_{na} \si g\}_{m,n\in \mz},$ so often it is a better strategy  to
construct  dual pairs
$\{E_{mb}T_{na}g\}_{m,n\in \mz}$,$ \{E_{mb}T_{na}h\}_{m,n\in \mz}$ such that $g$ and $h$
have required properties.
Dual pairs of Gabor frames have been characterized by
Ron \& Shen \cite{RoSh1} and  Janssen \cite{J}:

\bt  Two Bessel sequences  $\{E_{mb}T_{na}g\}_{m,n\in
\mz}$ and  $\{E_{mb}T_{na}h\}_{m,n\in \mz}$ form dual frames for $\ltr$ if
and only if \bes \sum_{k\in \mz} \overline{g(x-n/b-ka)}h(x-ka) =
b\delta_{n,0}, \ a.e. \ x\in [0,a].\ens \et

One of the most important results in Gabor analysis is the so-called
{\it duality principle.} It was discovered almost simultaneously by three
groups of researchers, namely  Daubechies \&
Landau \& Landau \cite{DLL}, Janssen \cite{Jan5}, and Ron \& Shen
\cite{RoSh1}. It concerns  the relationship between frame
properties for a function $g$ with respect to the {\it lattice}
\index{lattice}
$\{(na,mb)\}_{m,n\in \mz}$ and with respect to the so-called {\it
dual lattice} \index{dual lattice} $\{(n/b,m/a)\}_{m,n\in \mz}$:

\bt \label{223f} Given $g\in \ltr$ and $a,b>0,$  the following are equivalent:

 \bei \item[(i)]  $\mts$ is a frame for $\ltr$ with bounds
$A,B;$
\item[(ii)]  $\{\frac1{\sqrt{ab}} \,
E_{m/a}T_{n/b}g\}_{m,n\in \mz}$ is a Riesz sequence with bounds
$A, B$. \eni\et

The intuition behind the duality principle is that if $\mts$ is a frame
for $\ltr,$ then $ab\le 1,$ i.e., the sampling
points $\{(na,mb)\}_{m,n\in \mz}$ are ``sufficiently dense." Therefore
the points $\{(n/b,m/a)\}_{m,n\in \mz}$ are ``sparse," in the sense that
$\frac1{ab}\ge 1.$ Technically, this implies that the functions
$\{\frac1{\sqrt{ab}} \,
E_{m/a}T_{n/b}g\}_{m,n\in \mz}$ are linearly independent and only span a
subspace of $\ltr.$ The reason for the importance of the duality principle is that in general it is much easier to check that a system of vectors is a Riesz sequence than to check that it is a frame.
The duality principle is clearly related with the {\it Wexler-Raz theorem} stated next,
\index{Wexler-Raz thorem} which was discovered in 1994.

\bt
If the Gabor systems $\mts$ and
$\{E_{mb}T_{na}h\}_{m,n\in\mz}$ are Bessel sequences, then the following are
equivalent:

\bei \item[(i)] The Gabor systems $\mts$ and
$\{E_{mb}T_{na}h\}_{m,n\in\mz}$ are dual frames;
\item[(ii)] The Gabor
systems $\{\frac1{\sqrt{ab}} \, E_{m/a}T_{n/b}g\}_{m,n\in \mz}$
and $\{\frac1{\sqrt{ab}} \, E_{m/a}T_{n/b}h\}_{m,n\in \mz}$ are
biorthogonal, i.e.,
\bes \la \frac1{\sqrt{ab}} \, E_{m/a}T_{n/b}g,
\frac1{\sqrt{ab}} \, E_{m^\prime/a}T_{n^\prime/b}h\ra =
\delta_{m,m^\prime} \delta_{n, n^\prime}. \ens
\eni \et

\section{ An extension problem for wavelet frames} \label{404a}

\subsection{Background on the extension problem}

Extension problems have a long history in frame theory. It has been shown by several authors (see, e.g., \cite{CL} and \cite{LS}) that for any  Bessel sequence $\ftk$ in a separable Hilbert space $\h,$ there exists a sequence $\gtk$ such that
$ \ftk \, \cup \gtk$ is a tight frame for $\h.$ A natural generalisation to
construction of dual frame pairs appeared in \cite{CKK1}; we need to refer to the proof later, so we include it here as well.

\bt \label{223c} Let $\fti$ and $\gti$ be Bessel sequences in $\h$. Then there exist Bessel sequences $\ptj$ and $\qtj$ in $\h$ such that $\fti \cup \ptj$ and $\gti \cup \qtj$ form a pair of dual frames for $\h$. \et

\bp Let $T$ and $U$ denote the preframe operators  for $\fti$ and $\gti,$ respectively, i.e.,
\bes T,U: \ell^2(I)\to \h, \ T \{c_i\}_{i\in I}=\sui c_if_i, \ \ U \{c_i\}_{i\in I}=\sui c_ig_i.\ens Let $\atj, \btj$ denote any pair of dual frames for $\h.$ Then
\bes f= UT^*f+(I-UT^*)f & = & \sui \la f,f_i \ra g_i+
\suj \la  (I-UT^*)f, a_j \ra b_j \\ & = & \sui \la f,f_i \ra g_i+
\suj \la  f, (I-UT^*)^* a_j \ra b_j
\ens The sequences $\fti, \gti,$ and $\btj$ are Bessel sequences by definition, and one can verify that $\{(I-UT^*)^*a_j\}_{j\in J} $ is a Bessel sequence as well. The result now follows from Lemma \ref{223a}.\ep

The reason for the interest in this more general version of the frame extension is that it often is possible to construct dual pairs of frames with properties that are impossible for tight frames.
For example, Li and Sun showed in \cite{LS} that if
$ab\le 1$ and  $\mtgo$ is a Bessel sequences in  $\ltr,$  then there exists a Gabor systems $\mtgt$  such that $\mtgo \cup  \mtgt$ is a tight frame for $\ltr$. However, if we ask for extra properties of the functions $g_1$ and $g_2$ such an extension might be impossible. For example, if the given function $g_1$ has compact support, it is natural to ask for the function $g_2$ having compact support as well, but by \cite{LS} the existence of such a function is only guaranteed if  $| \supp g_1 | \le b^{-1}.$  On the other hand, such an extension can always be obtained in the setting of dual frame pairs \cite{CKK1}:

\bt
Let $\mtgo$ and $\mtho$ be Bessel sequences in  $\ltr,$ and
assume that $ab\le 1.$ Then the following hold:

\bei \item[(i)] There exist Gabor systems $\mtgt$ and $\mtht$
in $\ltr$ such that \bes \mtgo \cup \mtgt \ \mbox{and} \ \mtho \cup \mtht \ens form a
pair of dual frames for $\ltr$.
\item[(ii)] If $g_1$ and $h_1$ have compact support, the functions
$g_2$ and $h_2$ can be chosen to have compact support.
\eni \et

\bp Let us give the proof of (i). Let $T$ and $U$ denote the preframe operators for $\mtgo$ and $\mtho,$ respectively. Then
\bes UT^*f= \sumn \la f,E_{mb}T_{na}g_1 \ra E_{mb}T_{na}h_1.\ens
Consider the operator $\Phi:=I-UT^*,$  and let $\mtro,$ \ $\mtrt$ denote any pair of dual frames for $\ltr.$  By the proof of Theorem \ref{223c},  $\mtgo \cup \{\Phi^*E_{mb}T_{na}r_1\}_{m,n\in \mz} $ and
$\mtho \cup \{E_{mb}T_{na}r_2\}_{m,n\in \mz} $ are dual frames
for $\ltr.$ By Lemma \ref{223b} we know that $\Phi^*$ commutes with the time-frequency shift operators $E_{mb}T_{na}.$ This concludes the proof. \ep

\subsection{The extension problem for wavelet frames}

It turns out that the extension problem for wavelet systems is considerably more involved than for Gabor systems.
In order to explain this, consider the proof of Theorem \ref{223c}
and assume that $\fti$ and $\gti$ have wavelet structure, i.e.,
$\fti= \wpo$ and $\gti=\wpto$ for some $\psi_1, \widetilde{\psi_1}\in \ltr.$ Assume further that these sequences are Bessel sequences, with preframe operators $T,U,$ respectively. Then, still referring to the proof of Theorem \ref{223c},
$(I-UT^*)^*a_j= (I-TU^*)a_j.$
Unfortunately the operator $TU^*$ in general does not commute with
$D^jT_k,$ so even if we choose $\atj$ to have wavelet structure, the system $\{(I-TU^*)a_j\}_{j\in J}$ might not be a wavelet system.
Thus, we can not apply the proof technique from the Gabor case.
The following partial result was obtained in \cite{CKK1}.

 \bt \label{213a} Let $\wpo$ and $\wpto$ be Bessel sequences in
$\ltr$. Assume that the Fourier transform of $\widetilde{ \psi_1}$
satisfies
    \begin{equation} \label{223d}
        \supp\  \widehat {\widetilde{\psi_1}}\subseteq [-1,1].
    \end{equation}
Then there exist wavelet systems $\wpt$ and $\wptt$ such that \bes
\wpo \cup \wpt \ \mbox{and} \ \wpto \cup \wptt \ens form dual
frames for $\ltr.$ If we further assume that
$\widehat{\psi_1}$ is compactly supported and that
    \begin{equation*}
        \supp\  \widehat {\widetilde{\psi_1}}\subseteq [-1,1]\setminus [-\epsilon, \epsilon]
    \end{equation*} for some $\epsilon>0, $ the functions
$\psi_2 $ and $ {\widetilde{\psi_2}}$ can be chosen to
have compactly supported Fourier transforms  as well. \et

In the Gabor case, no assumption of compact support was necessary, neither for the given functions nor their Fourier transform. From this point of view
it is natural to ask whether the assumption \eqref{223d} is necessary in
Theorem \ref{213a}.

\vspace{.1in}\noindent{\bf Questions:}  {\it Let $\wpo$ and $\wpto$ be Bessel sequences in
$\ltr$.
\bei \item[(i)]
Does there exist functions $\psi_2, \widetilde{\psi_2}\in \ltr$  such that \bee \label{287a}
\wpo \cup \wpt \ \mbox{and} \ \wpto \cup \wptt \ene form dual
frames for $\ltr?$
\item[(ii)] If $\widehat{\psi_1}$  and $ \widehat {\widetilde{\psi_1}}$ are compactly supported, can we find compactly supported functions
$\psi_2 $ and $\widetilde{\psi_2}\in \ltr$  such that the functions in \eqref{287a} form dual frames?
\eni}

\vspace{.1in} The problem (i) can also be formulated in the negative way:
can we find just one example of a pair of
Bessel sequences $\wpo$ and $\wpto$ that can not be extended to a pair of
dual wavelet frames, each with 2 generators?
The open question is strongly connected to the following conjecture by Han
\cite{Han}:

\vspace{.1in}\noindent{\bf Conjecture by Deguang Han:}
{\it Let $\{D^jT_k\psi_1\}_{j,k\in \mz}$
be a wavelet frame with upper frame bound $B.$ Then there exists $D>B$ such
that for each $K\ge D,$ there exists $\widetilde{\psi_1}\in \ltr$ such that
$\{D^jT_k\psi_1\}_{j,k\in \mz} \cup \{D^jT_k\widetilde{\psi_1}\}_{j,k\in \mz}$
is a tight frame for $\ltr$ with bound $K.$}

\vspace{.1in} The paper \cite{Han} contains an example showing that (again in contrast with the Gabor setting) it might not be possible to extend
the Bessel system $\{D^jT_k\psi_1\}_{j,k\in \mz}$ to a tight frame
without enlarging the upper bound; hence it is essential that the conjecture
includes the option that the extended wavelet system has a strictly larger frame bound than the upper frame bound $B$ for $\{D^jT_k\psi_1\}_{j,k\in \mz}.$
We also note that Han's conjecture is based on an example, where
$ \supp\  \widehat {\psi_1}\subseteq [-1,1],$ i.e., a case that is covered by Theorem \ref{213a}.

Observe that a pair of wavelet Bessel sequences always can be extended to dual wavelet frame pairs by adding {\it two} pairs of wavelet systems. In fact, we can always add one pair of wavelet systems
that cancels the action of the given wavelet system, and another one that yields a dual pair of wavelet frames by itself. Thus, the issue is really whether it is enough to add one pair of wavelet systems, as stated in the formulation of the open problem.

Note that extension problems have a long history in frame theory. Most of the results are dealing with the unitary extension principle \cite{RoSh2},
\cite{RoSh3}
and its variants, and are thus based on the assumption of an underlying refinable function. The open problems formulated in this section are not based on such an assumption.

\section{ The duality principle in Hilbert spaces} \label{404c}

\subsection{The background}

The duality principle, Theorem \ref{223f}, is one of the key results in Gabor analysis. Therefore it is natural to ask whether a similar result is valid for general frames in Hilbert spaces. Casazza,  Kutyniok, and Lammers investigated this in \cite{CKL}, and
introduced the  {\it R-dual} of a sequence $\fti$ in a general separable Hilbert space $\h:$

\bd  Let $\eti$ and $\hti$ denote
orthonormal bases for $\h,$ and let $\fti$ be any sequence in $\h$
for which $\sui |\la f_i,e_j\ra|^2 < \infty$ for all $j\in I.$ The
R-dual of $\fti$ with respect to the orthonormal bases $\eti$ and
$\hti$ is the sequence $\otj$ given by \bee \label{207a} \oj= \sui
\la f_i, e_j\ra h_i, \ i\in I.\ene \ed

Let us state some of the central results from the paper \cite{CKL}.

\bt \label{213b} Define the R-dual $\otj$ of a
sequence $\fti$ as above. Then the following hold:

\bei \item[(i)] For all $i\in I,$ \bee \label{207d} f_i= \suj \la
\oj, h_i\ra e_j,\ene i.e., $\fti$ is the R-dual sequence of $\otj$
w.r.t. the orthonormal bases $\hti$ and $\eti.$
\item[(ii)] $\fti$ is a Bessel sequence with bound $B$ if and only $\oti$ is a
Bessel sequence with bound $B.$
\item[(iii)] $\fti$ satisfies the lower frame condition with bound
$A$ if and only if $\otj$ satisfies the lower Riesz sequence
condition with bound $A$.
\item[(iv)] $\fti$ is a frame for $\h$ with bounds $A,B$ if and
only if $\otj$ is a Riesz sequence in $\h$ with bounds $A,B.$
\item[(v)] Two Bessel sequences $\fti$ and $\gti$ in $\h$ are
dual frames if and only if the associated R-dual sequences $\otj$
and $\{ \gamma_j\}_{j\in I}$  with respect to the same orthonormal bases $\eti$ and $\hti$ satisfy that \bee
\label{207f} \la \oj, \gamma_k\ra = \delta_{j,k}, \ j,k\in I.\ene
\eni \et

From Theorem \ref{213b} it is evident that the relations between a given sequence $\fti$ and its R-dual $\otj$ resembles the relations between
a Gabor system $\mts$ and the corresponding Gabor system
$\{\frac1{\sqrt{ab}} \, E_{m/a}T_{n/b}g\}_{m,n\in
\mz}$ on the dual lattice, see Theorem \ref{223f}.  This raises the natural question whether the duality principle is actually a special case of Theorem \ref{213b}.  That is, can $\{\frac1{\sqrt{ab}} \, E_{m/a}T_{n/b}g\}_{m,n\in
\mz}$ be realized as the R-dual of $\mts$ with respect to certain choices
of orthonormal bases $\{e_{m,n}\}_{m,n\in \mz}$ and $\{h_{m,n}\}_{m,n\in \mz}$?

The paper \cite{CKL} does not provide the complete answer to this question, but it was clearly the driven force of the authors.
The paper contains the following partial results:

\bei \item If $\mts$ is a frame and
$ab=1,$ then $\{\frac1{\sqrt{ab}} \,
E_{m/a}T_{n/b}g\}_{m,n\in \mz}$ can be realized as the R-dual of
$\mts$ with respect to certain choices of orthonormal bases $\{e_{m,n}\}_{m,n\in \mz}$ and
$ \{h_{m,n}\}_{m,n\in \mz} $.
\item If $\mts$ is a
tight frame, then $\{\frac1{\sqrt{ab}} \,
E_{m/a}T_{n/b}g\}_{m,n\in \mz}$ can be realized as the R-dual of
$\mts$ with respect to certain choices of orthonormal bases $\{e_{m,n}\}_{m,n\in \mz}$ and
$ \{h_{m,n}\}_{m,n\in \mz} $.
\eni

We note that a complementary
approach to the question was given \cite{CKK2}, where the authors ask for general  conditions on two
sequences $\fti, \otj$ such that $\otj$ is the R-dual of $\fti$
with respect to {\it some} choice of the orthonormal bases $\eti$
and $\hti.$
The following result is proved in \cite{CKK2}.

\bt  \label{223g} Let $\otj$ be a
Riesz sequence spanning a proper subspace $W$ of $\h$ and $\eti$
an orthonormal basis for $\h$. Given any frame $\fti$ for $\h$
the following are equivalent:
\bei \item[(i)] $\otj$ is an R-dual of $\fti$  with respect to $\eti$ and
some orthonormal basis $\hti$.
\item[(ii)] There exists an orthonormal basis $\hti$ for $\h$
satisfying
\bee \label{197a} f_i= \suj \la \oj, h_i\ra e_j, \ \forall i\in I.\ene
\item[(iii)] The sequence $\nti$ given by
\bee \label{104a}
n_i:= \suk \la e_k, f_i\ra \ottk, \ i\in I.\ene
is a
tight frame for $W$ with frame bound $E=1.$ \eni \et

In \cite{CKK2} it is shown that in the setting of Theorem \ref{223g},
the sequences $\hti$ satisfying \eqref{197a} are
characterized as $h_i=m_i+n_i,$ where $m_i\in \wp.$
Also,  if $\otj$ has the bounds $C,D$ and $\fti$ has the bounds $A,B,$ then $\nti$ is always a frame for $W,$ with bounds $A/D, B/C;$ the only question is whether it can be made tight by an appropriate choice of the orthonormal basis $\eti.$

We can formulate the main questions as follows:

\vspace{.1in}\noindent
{\bf Questions:} {\it
\bei \item[(i)] Can the duality principle in Gabor analysis be realized via the theory for R-duals? That is, given any Gabor frame $\mts,$ can the
Riesz sequence  $\{\frac1{\sqrt{ab}} \,
E_{m/a}T_{n/b}g\}_{m,n\in \mz}$  be realized as the R-dual of $\mts$?
\item[(ii)] What are the general conditions on $\fti$ and $\otj$ such that there is an orthonormal basis  $\eti$ with the property that $\nti$ in \eqref{104a} is a tight frame for $W$ with bound $A=1$?
    \item[(iii)] Does the connections between the duality principle and the R-duals lead to useful results for other structured systems, e.g., wavelet systems?
\eni}

\section{Wavelet packet frames} \label{404d}
\subsection{The background}

Let
$\{a_j\}_{j\in \mz}$ be a collection of positive numbers,
let $b>0$, and let $\{c_m\}_{m\in \mz}$ be a collection
of points in $\mr$. Given a function $g\in \ltr,$ we will consider
the system of functions
\bee \label{w} \{D_{a_j}T_{kb}E_{c_m}g\}_{
j,m,k\in \mz}.\ene
A system of the type \eqref{w} is called a {\it wave packet system}.
Note that
\bes D_{a_{j}}T_{bk}E_{c_{m}}g=T_{a_{j}^{-1}kb}D_{a_{j}}E_{c_{m}}g.\ens
Thus,  a wave packet system is a special case of a generalized shift-invariant system, as considered, e.g., in \cite{HLWI} and \cite{RS6}. We note that wave packet systems have also been considered by Czaja, Kutyniok and Speegle in \cite{CzK}. They
proved that certain geometric
conditions on the set of parameters in a wave packet systems are
necessary in order for the system to form a frame, and also provide constructions of frames and orthonormal bases, based on characteristic functions.

We first state a result from \cite{CR} about construction of frames of the form $\{D_{a_j}T_{kb}E_{c_m}g\}_{
j,m,k\in \mz},$ .

\bt \label{t1} Let $\{a_{j}\}_{j\in\mz}$ be a sequence of positive
numbers , $b>0$, $\{c_{m}\}_{m\in\mz}\subset \mr$ and $g\in\ltr$.
Assume that
\begin{equation}\label{bzm1}
B:=\frac{1}{b}\sup_{\gamma\in\mr}\sum_{j,m\in\mz}\sum_{k\in\mz}|\hat{g}(a_j^{-1}\gamma
-c_{m})\hat{g}(a_j^{-1}\gamma-c_{m}- k/b)|<\infty.
\end{equation}
Then the wave packet system
$\{D_{a_{j}}T_{bk}E_{c_{m}}g\}_{j,m\in\mz,k\in\mz}$ is a Bessel
sequence with bound $B$. Further, if also
\begin{eqnarray*}
A&:=&\frac{1}{b}\inf_{\gamma\in\mr}\bigg(\sum_{j,m\in\mz}|
\hat g (a_j^{-1}\gamma-c_{m})|^{2}\\
&-& \sum_{0\neq k\in \mz}\sum_{j,m\in\mz}|\hat g
(a_j^{-1}\gamma-c_{m})\hat g (a_j^{-1}\gamma-c_{m}- k/b)|\bigg)>0,
\end{eqnarray*} then the wave packet system
$\{D_{a_{j}}T_{bk}E_{c_{m}}g\}_{j,m\in\mz,k\in\mz}$ is a frame for
$\ltrd$ with bounds $A$ and $B$. \et

The Bessel condition puts certain restrictions on the
numbers $\{a_j \}_{j\in \mz} $ and the
distribution of the points $\{c_m\}_{m\in \mz}$.
Let us consider a function $g\in\ltr$ such
that for some interval $I\subset \mr$ we have \bee \label{g}
|\hat{g}(\gamma)|>\epsilon>0, \ \ga \in I. \ene Also, we will assume that the sequence
of points
$\{c_{m}\}_{m\in\mz}\subset \mr$ is such that for some $r>0$ \bee \label{c} \bigcup_{m\in\mz}(c_{m}+[0,r])=\mr.\ene For many natural choices of the numbers $a_j$, it was shown in \cite{CR} that the above
conditions on $g\in \ltr$ and $\{c_{m}\}_{m\in\mz}$ exclude the
frame property for the associated wave packet system. Similar results appeared in
\cite{CzK}, formulated in terms of the upper Beurling density.

\bl \label{cr} Let $ \{a_j\}_{j\in \mz}$ be a sequence of positive numbers,
$b>0$ and
assume that there exists a number $C>0$ and an infinite index set
$J\subseteq \mz$ such that \bee \label{ff} a_j\le C, \
\forall j\in I. \ene Assume that $\{c_{m}\}_{m\in\mz}\subset
\mr$ satisfies \eqref{c}. Then no function $g\in\ltr$
satisfying \eqref{g} can generate a Bessel sequence
$\{D_{a_{j}}T_{kb}E_{c_{m}}g\}_{j,m,k\in\mz}$. \el

As a special case, the assumptions on $g$ and $\{c_{m}\}_{m\in\mz}$
 in Lemma \ref{cr} exclude the frame property for the system
$\{D_{a^{j}}T_{kb}E_{c_{m}}g\}_{j,m,k\in\mz}$  for any $a>1.$ In other words, if we
want
$\{D_{a^{j}}T_{kb}E_{c_{m}}g\}_{j,m,k\in\mzd}$ to be a frame under the weak assumption \eqref{g}, we
need to consider sparser distributed points $\{c_m\}_{m\in \mz}$
than the one considered in \eqref{c}. But it is not known how to identify suitable conditions of the function $g$ and the distribution of the points $\{c_m\}_{m\in \mz}$ to ensure the frame property.

Let us continue to discuss  wave packet systems of the special form  \\
$\{D_{a^{j}}T_{bk}E_{c_{m}}\psi\}_{j\in\mz,k\in\mz, m\in \mz}.$
 Let \bes \Lambda:= \left\{\frac{a^jn}{b} \ \big|
j,n\in \mz\right\};\ens and, given $\alpha\in \Lambda,$ let \bes
J_\alpha:= \{j\in \mz \big| \ \exists n\in \mz \ \mbox{such that}
\ \alpha= \frac{a^jn}{b}\}. \ens Finally,
let \bes {\cal D}= \{f\in \ltr \ \big| \ \hat{f}\in L^\infty(\mr),
\ \supp \, \hat{f} \ \mbox{is compact}\}.\ens
Then the
main result in \cite{HLWI} takes the following form:

\bt \label{104b} Let $\psi, \tpi\in \ltr,$ and assume that
$\{D_{a^{j}}T_{bk}E_{c_{m}}\psi\}_{j\in\mz,k\in\mz, m\in \mz}$ and
$\{D_{a^{j}}T_{bk}E_{c_{m}}\tpi\}_{j\in\mz,k\in\mz, m\in \mz}$ are
Bessel sequences. If \bee \label{lic} L(f):= \sum_{j\in \mz}
\sum_{m\in \mz} \sum_{n\in \mz} \int_{\supp \, \hat{f}}
|\hat{f}(\ga + \frac{a^jn}{b})|^2  \, | \hpi(a^{-j}\ga-c_m)|^2 \,
d\ga <\infty \ene for all $f\in {\cal D},$ then
$\{D_{a^{j}}T_{bk}E_{c_{m}}\psi\}_{j\in\mz,k\in\mz, m\in \mz}$ and
$\{D_{a^{j}}T_{bk}E_{c_{m}}\tpi\}_{j\in\mz,k\in\mz, m\in \mz}$
form pairs of dual frames for $\ltr$ if and only if \bee
\label{g1} \sum_{j\in J_\alpha} \sum_{m\in \mz}
\hat{\psi}(a^{-j}\ga -c_m) \overline{\htpi(a^{-j}(\ga+\alpha)
-c_m)} = b\, \delta_{\alpha,0}.\ene \et

Theorem \ref{104b} leads to a sufficient condition for duality for two wave
packet systems that resembles the  versions we have seen for Gabor systems and wavelet systems:

\bc \label{104d} Assume that $\{D_{a^{j}}T_{bk}E_{c_{m}}\psi\}_{j\in\mz,k\in\mz, m\in
\mz}$ and $\{D_{a^{j}}T_{bk}E_{c_{m}}\tpi\}_{j\in\mz,k\in\mz, m\in
\mz}$ are Bessel sequences and that the local integrability condition
\eqref{lic} holds. Then $\{D_{a^{j}}T_{bk}E_{c_{m}}\psi\}_{j\in\mz,k\in\mz, m\in
\mz}$ and $\{D_{a^{j}}T_{bk}E_{c_{m}}\tpi\}_{j\in\mz,k\in\mz, m\in
\mz}$ are dual frames if the following two conditions hold:

\bee \label{c1} \sum_{j\in \mz} \sum_{m\in \mz} \hat{\psi}(a^{-j}\ga -c_m)
\overline{\htpi(a^{-j}\ga -c_m)} = b, \ a.e. \ \ga \in \mr, \\
\label{c2} \hat{\psi}(\ga) \overline{\htpi(\ga+ q)} = 0, \ a.e. \
\ga \in \mr \ \mbox{for} \ q\in b^{-1}(\mz\setminus \{0\}).\ene
\ec \bp The condition \eqref{c1} corresponds to the condition
\eqref{g1} with $\alpha= 0.$ For $\alpha \neq 0,$ we note that
\bes j\in J_\alpha \Leftrightarrow \exists n\in \mz: \alpha =
\frac{a^j}{b}n \Leftrightarrow \alpha a^{-j} \in b^{-1}\mz.\ens
Thus, the condition \eqref{c2} implies that the double sum in
\eqref{g1} vanishes for all $j\in J_\alpha$. \ep

Similar to our discussion of the frame property, it is not known how to identify suitable conditions of the function $g$ and the distribution of the points $\{c_m\}_{m\in \mz}$ to ensure the duality property.  Let us collect the raised
problems:

\vspace{.1in}
\noindent{\bf Problem:}

{\it

\bei \item[(i)] Find directly applicable conditions on a function $g$ and the distribution of the sequences  $ \{a_j\}_{j\in \mz}$ and
$\{c_{m}\}_{m\in\mz}$ such that $\{D_{a_{j}}T_{kb}E_{c_{m}}g\}_{j,m,k\in\mz}$
is a frame for some $b>0.$
\item[(ii)] Find directly applicable conditions on a function $g$ and the distribution of the sequence
$\{c_{m}\}_{m\in\mz}$ such that $\{D_{a_{j}}T_{kb}E_{c_{m}}g\}_{j,m,k\in\mz}$
is a frame for some $a,b>0.$
\item[(iii)]Find possible distributions of numbers $ \{a_j\}_{j\in \mz}$  and
$\{c_{m}\}_{m\in\mz}\subset \mr$ such that dual frames of the form
$\{D_{a_{j}}T_{kb}E_{c_{m}}g\}_{j,m,k\in\mz}$ and
$\{D_{a_{j}}T_{kb}E_{c_{m}}h\}_{j,m,k\in\mz}$ exist.
\eni}

As mentioned, the paper \cite{CzK} contains some answers to the raised questions in the special case where $g$ is chosen to be a characteristic function. But the are no results available for more general functions.

\section{B-splines and Gabor frame} \label{404e}
We will now consider frame properties for Gabor systems generated by the B-splines, defined inductively for $N\in \mn$ by \bes B_1(x)=
\chi_{[0,1]}(x), \ \ B_{N+1}(x)=B_N*B_1(x)= \int_0^1 B_N(x-t)\, dt.\ens

The B-splines are well studied, and have many desirable properties. For example,  \bei
\item[(i)] $supp \ B_N = [0, N]$ and
$B_N>0$ on $]0, N[$.
\item[(ii)] $\inr B_N(x)dx=1$.
\item[(iii)] $ \sum_{k\in \mz}B_N(x-k)=1$
\item[(iv)] For any $N\in \mn$,
$\widehat{B_N}(\ga)=\left(\frac{1- e^{-2\pi i \ga}}{2\pi i
\ga}\right)^N.$ \eni

The connection of the B-splines and modern harmonic analysis is well known. For example, the B-splines lead to constructions of tight wavelet frames via the unitary extension principle by Ron and Shen \cite{RoSh3}.
It is a classical result that for  $N\in \mn$, the B-spline $B_N$ generates a Gabor
frame $\{E_{mb}T_nB_N\}_{m,n\in \mz}$ for all $(a,b)\in ]0, N[\times ]0, 1/N[$. It is also known that for
$b\le \frac1{2N-1}$ and $a=1,$ the Gabor frame $\{E_{mb}T_nB_N\}_{m,n\in \mz}$
has a dual frame $\{E_{mb}T_nh_N\}_{m,n\in \mz},$ for which the function $h_N$ is just a (finite) linear combination of shifts of $B_N,$ see \cite{C19, CR}.

However, the exact range
of parameters $(a,b)$ for which $B_N$ generates a frame is unknown.

\vspace{.2in}\noindent{\bf Question:} {\it Given $N\in \mn,$ characterize  the $(a,b)\in
\mr^2$ for which $B_N$ generates a Gabor frame $\{E_{mb}T_{na}B_N\}_{m,n\in \mz}$.}

\vspace{.2in}
A characterization of the parameters $(a,b)\in \mr$ that yield Gabor frames
$\mts$ is only known for a few types of functions $g,$ including the Gaussian \cite{Ly}, \cite{Se2}, \cite{SW}, the hyperbolic secant \cite{JS2}, and the class of totally positive functions of finite type \cite{GS}.

The exact answer to the above question is bound to be complicated. For example, consider the
B-spline $B_2.$ It is easy to show that   $\{E_{mb}T_{na}B_2\}_{m,n\in \mz}$ can not be a
frame for any $b>0$ whenever $a\ge 2$.  On the other hand, it was shown
in \cite{GJKP} that for $b=2,3,\dots,$
$\{E_{mb}T_{na}B_2\}_{m,n\in \mz}$ can not be a frame for any $a>0$.

We note that it is the lower frame condition that causes the problem. In fact, by Lemma \ref{804a} we know that
$\{E_{mb}T_{na}B_N\}_{m,n\in \mz}$ is a Bessel sequence for all $a,b>0$ and all $N\in \mn.$

The difficulty of the problem is illustrated by the
related problem of characterizing all $a,b,c>0$ for which
$\{E_{mb}T_{na}\chi_{[0,c]}\}_{m,n\in \mz}$ is a frame for $\ltr.$
Janssen considered this problem in the paper \cite{Jan7}, and gave the answer in 8 special cases.
The full problem of characterizing $a,b,c>0$ for which
$\{E_{mb}T_{na}\chi_{[0,c]}\}_{m,n\in \mz}$ is a frame was finally solved by  Dai and Sun, see \cite{DS}.

\section{Finite structured frames} \label{404g}

Even for frames in infinite-dimensional Hilbert spaces like $\ltr,$ concrete implementations always has to take place on finite subcollections. However, it is well known that finite subfamilies might have properties that are quite different from the full system. For example, any finite collection
of vectors $\{f_k\}_{k=1}^N$ in a Hilbert space is a frame for its linear span, while
an infinite collection of vectors $\ftk$ certainly does not need to form a frame for
$\span \ftk.$

 As  motivation for the following, let us state a result by Kim and Lim \cite{KL}, see also \cite{C16}.

\bl \label{804b} Assume that $\ftk$ is an overcomplete frame, for which any finite
subfamily $\{f_k\}_{k=1}^N$ is linearly linearly independent. Let $A_N$  denote a lower bound for  $\{f_k\}_{k=1}^N,$ considered as a frame for
the space $\h_N:= \mbox{span} \{f_k\}_{k=1}^N.$ Then
\bes A_N\to 0 \ \mbox{as} \ N\to \infty.\ens \el
Lemma \ref{804b} has an interesting consequence for Gabor frames
$\mts$ with $ab<1:$ in fact, for such systems, it will follow from our discussion in Section \ref{804c} that the lower frame bound for any subset
$\{E_{mb}T_{na}g\}_{|m|, |n| \le N}$  tend to zero as $N\to \infty.$

We will now discuss two questions that arise naturally when considering Gabor systems, respectively, collections of complex exponentials.

\subsection{The Heil--Ramanathan--Topiwala conjecture} \label{804c}
So far, we have only discussed what could be called {\it regular} Gabor systems, meaning that the time-frequency shifts form a lattice
$\{(na,mb)\}_{m,n\in \mz}$ in $\mr^2.$ I is also possible (though considerably more complicated) to consider Gabor systems with arbitrary time-frequency shift, i.e., systems of the form
$\{e^{2\pi i \lambda_nx}g(x-\mu_n)\}_{n\in I}$ for some collection of points
$\{(\lambda_n, \mu_n)\}_{n\in I}.$ An interesting and surprisingly difficult problem for such systems was formulated by Heil, Ramanathan and Topiwala in the paper \cite{HRT} from 1996:

\vspace{.1in}\noindent{\bf Conjecture by Heil, Ramanathan and Topiwala:} {\it Let $\{(\lambda_n, \mu_n)\}_{n\in I}$ denote a finite set of distinct points in $\mr^2,$ and let $g\in \ltr$ be an arbitrary nonzero function. Then the finite Gabor system $\{e^{2\pi i \lambda_nx}g(x-\mu_n)\}_{n\in I}$ is linearly independent.}

\vspace{.1in} The conjecture  is still open. So far, it has been confirmed in several important special cases. For example, Linnell \cite{Linn} proved the conjecture in the case where
$\{(\lambda_n, \mu_n)\}_{n\in I}=\{(mb,na)\}_{n=1, m=1}^{N,M},$ i.e., for subsets of the lattice-type Gabor systems described in Section \ref{273a}.

A detailed account of the conjecture and the various known results can be found in the paper \cite{H} by Heil. As inspiration for the reader, we just mention one seemingly easy case, taken from \cite{H}, where the conjecture is still open:

\vspace{.1in}\noindent{\bf Special case of the conjecture:} {\it Let $g\in \ltr$ be a nonzero function. Then the set of functions
\bes \{g(x), g(x-1), e^{2\pi i x} g(x), e^{2\pi i \sqrt{2}x}g(x-\sqrt{2})\}\ens is linearly independent.}

\subsection{Lower bounds for a finite collection of  exponentials}
Motivated by the fact that the (scaled) exponentials $\{\frac1{\sqrt{2\pi}}\, e^{ inx}\}_{n\in \mz}$ forms an orthonormal basis for $\ltp,$ it is natural to consider {\it nonharmonic Fourier series,} i.e., expansions in terms of exponentials
$\{e^{i\lambda_nx}\}_{n\in \mz}$
for a collection of numbers $\{\lambda_n\}_{n\in \mz}.$ Analysis of such systems is actually a classical subject, and we refer to the excellent presentation in \cite{Y}.

Consider now a finite  collection of distinct real number $\{\lambda_n\}_{n=1}^N$ ordered increasingly,
\bes \lambda_1 < \lambda_2 < \cdots <\lambda_N.\ens Then $\{e^{i\lambda_nx}\}_{n=1}^N$ is a frame for its linear span in $\ltp.$
Letting $A_N$ denote a lower frame bound
for the frame sequence $\{e^{i\lambda_nx}\}_{n=1}^N,$ it follows from Lemma \ref{804b} that $A_N\to 0$ as
$N\to \infty.$ The decay of $A_N$ is estimated in \cite{CL1}, where
it is showed that if we choose  $\delta \le 1$ such that $|\lambda_k-\lambda_{k+1}|
\ge \delta$ for all $k=1, \dots, N-1,$ then
\bes A_N:= 1.6 \cdot 10^{-14} \left( \frac{\delta}{2} \right)^{2N+1} \frac1{((N+1)!)^8}\ens is a lower frame bound for $\{e^{i\lambda_nx}\}_{n=1}^N.$ This is obviously a very crude estimate, which leads to a natural question:

\vspace{.1in}\noindent{\bf Open question:} {\it How can better estimates for the lower frame bound of $\{e^{i\lambda_nx}\}_{n=1}^N$ be obtained?}

\section{The Feichtinger conjecture} \label{404h}

Around 2002 Feichtinger observed that several Gabor frames in the literature could  be split into finite collections of Riesz sequences. He
formulated the following conjecture in emails to some colleagues in harmonic analysis:

\vspace{.1in}\noindent{\bf The Feichtinger conjecture:} {\it Let $\ftk$ be a frame with the property that $\inf_{k\in \mn} ||f_k||>0.$ Then $\ftk$ can be partitioned into a finite union of Riesz sequences.}

\vspace{.1in} Relatively soon, the first positive partial results were published in \cite{G5} and \cite{CCLV}. However, the general question turned out to be very difficult. Around 2005 it was shown by Casazza and Tremain that the Feichtinger conjecture is equivalent to the  Kadison--Singer conjecture from 1959, in the sense that either both conjectures are true or both are false. Later, Casazza related the conjecture to several other open problems in the literature. We refer to \cite{CE} and \cite{Cas} for  detailed descriptions of these conjectures. Shortly before submission of the current manuscript, the Feichtinger conjecture was reported to be solved affirmatively by Marcus, Spielman and Srivastava, see \cite{MSS}.

{\bf \vspace{.1in}
\noindent Ole Christensen\\
Department of Mathematics and Computer Science\\
Technical University of Denmark\\
Building 303 \\
2800 Lyngby  \\
Denmark \\
 Email: ochr@dtu.dk

 }
\end{document}